\input amstex
\documentstyle{amsppt}
\magnification=\magstep1                        
\hsize6.5truein\vsize8.9truein                  
\NoRunningHeads
\loadeusm

\magnification=\magstep1                        
\hsize6.5truein\vsize8.9truein                  
\NoRunningHeads
\loadeusm

\document
\topmatter

\date June 9, 2014
\enddate

\title Inequalities for Lorentz Polynomials 
\endtitle

\rightheadtext{Inequalities for Lorentz Polynomials}

\author Tam\'as Erd\'elyi
\endauthor

\address Department of Mathematics, Texas A\&M University,
College Station, Texas 77843 (T. Erd\'elyi) \endaddress

\email terdelyi\@math.tamu.edu (T. Erd\'elyi)
\endemail

\email terdelyi\@math.tamu.edu
\endemail

\endthanks


\dedicatory
\enddedicatory

\abstract
We prove a few interesting inequalities for Lorentz polynomials. A highlight of this paper 
states that the Markov-type inequality  
$$\max_{x \in [-1,1]}{|f^\prime(x)|} \leq n \max_{x \in [-1,1]}{|f(x)|}$$
holds for all polynomials of degree at most $n$ with real coefficients for which $f^\prime$ 
has all its zeros outside the open unit disk. Equality holds only for 
$f(x) := c((x \pm 1)^n - 2^{n-1})$  with a constant $0 \neq c \in {\Bbb R}$. This should be 
compared with  Erd\H os's classical result stating that  
$$\max_{x \in [-1,1]}{|f^\prime(x)|} \leq \frac n2 \left( \frac{n}{n-1} \right)^{n-1} \max_{x \in [-1,1]}{|f(x)|}$$ 
for all polynomials of degree at most $n$ having all their zeros in ${\Bbb R} \setminus (-1,1)$.
\endabstract

\endtopmatter

\head 1. Introduction \endhead

Let ${\Cal P}_n$ denote the collection of all polynomials of degree at most $n$  
with {\it real coefficients}. Let ${\Cal P}_n^c$ denote the collection of all polynomials 
of degree at most $n$ with {\it complex coefficients}. Let 
$$\|f\|_A := \sup_{x \in A}{|f(x)|}$$ 
denote the supremum norm of a complex-valed function $f$ defined on a set $A$.
The Markov inequality asserts that
$$\|f^\prime\|_{[-1,1]} \leq n^2 \|f\|_{[-1,1]}$$
holds for all $f \in {\Cal P}_n^c$.
The inequality
$$|f^\prime(x)| \leq \frac{n}{\sqrt{1-x^2}} \, \|f\|_{[-1,1]}$$
holds for all $f \in {\Cal P}_n^c$ and for all $x \in (-1,1)$, 
and is known as Bernstein inequality. For proofs of these see [2] or [5], 
for instance.  
Various analogues of the above two inequalities are known in which the
underlying intervals, the maximum norms, and the family
of functions are replaced by more general sets, norms, and families of
functions, respectively. These inequalities are called Markov-type and
Bernstein-type inequalities.
If the norms are the same in both sides, the inequality is called
Markov-type, otherwise it is called Bernstein-type (this distinction
is not completely standard). Markov- and Bernstein-type inequalities
are known on various regions of the complex plane and the $n$-dimensional
Euclidean space, for various norms such as weighted $L_p$ norms, and
for many classes of functions such as polynomials with various constraints,
exponential sums of $n$ terms, just to mention a few. Markov- and
Bernstein-type inequalities have their own intrinsic interest. In addition,
they play a fundamental role in approximation theory.

It had been observed by Bernstein that Markov's
inequality for monotone polynomials is not essentially better than for
arbitrary polynomials. Bernstein proved that 
$$\sup_{f}{\frac{\|f^{\prime}\|_{[-1,1]}}{\|f\|_{[-1,1]}}} =
\left\{ \aligned \textstyle{\frac 14} (n+1)^2 \,, \qquad &\text{if} \enskip n \enskip \text{is odd} \\
\textstyle{\frac 14} n(n+2)\,, \qquad &\text{if} \enskip n \enskip \text{is even} \,, \endaligned \right. $$
where the supremum is taken over all $f \in {\Cal P}_n$ which are monotone on $[-1,1]$. 
See [22], for instance. 
This is surprising, since one would expect that if a polynomial is this far away 
from the ``equioscillating" property of the Chebyshev polynomial $T_n$, then 
there should be a more significant improvement in the Markov inequality. In [16] 
Erd\H os gave a class of restricted polynomials for which 
the Markov factor $n^2$ improves to $cn$. He proved that there is an absolute 
constant $c$ such that 
$$|f^{\prime}(x)| \leq \min \left\{\frac{c\,\sqrt{n}}{\left(1-x^2\right)^2}\,,
\enskip \frac{en}{2} \right\} \,\|f\|_{[-1,1]}\,, \qquad x \in (-1,1)\,,$$
for all $f \in {\Cal P}_n$ having all their zeros in
${\Bbb R} \setminus (-1,1)$.
This result motivated several people to study Markov- and Bernstein-type 
inequalities for polynomials with restricted zeros and under some other 
constraints. Generalizations of the above Markov- and Bernstein-type
inequality of Erd\H os have been extended in various directions by several people 
including Lorentz [20], Scheick [23], Szabados [24], M\'at\'e [21], 
P. Borwein [1], Erd\'elyi [6,7,9,12,13], Rahman and Schmeisser [22], 
Kro\'o and Szabados [18,19], Hal\'asz [17], and the list  can be even longer.  
A special attention is paid to the classes ${\Cal P}_{n,k}$ and ${\Cal P}_{n,k}^c$, 
where ${\Cal P}_{n,k}$ denotes the set of all polynomials of degree at most $n$ with 
{\it real coefficients} and with at most $k$ ($0 \leq k \leq n$) zeros in the open unit 
disk, and ${\Cal P}_{n,k}^c$ denotes the set of all polynomials of degree at most $n$ with 
{\it complex coefficients} and with at most $k$ ($0 \leq k \leq n$) zeros in the open unit 
disk. Associated with $0 \leq k \leq n$ and $x \in (-1,1)$, let
$$B_{n,k,x} := \sqrt{\frac{n(k+1)}{1-x^2}}\,, \qquad \quad
B_{n,k,x}^* := \max \left\{\sqrt{\frac{n(k+1)}{1-x^2}}\,, \enskip
n \log \left( \frac{e}{1-x^2} \right) \right\}\,,$$
and
$$M_{n,k} := n(k+1)\,, \qquad \quad M_{n,k}^* := \max\{n(k+1), \enskip n\log n\}\,.$$
In [10] and [11] it is shown that
$$c_1 \min\{B_{n,k,x}^*, M_{n,k}^*\}
\leq \sup_{f \in {\Cal P}_{n,k}^c}{\frac {|f^\prime(x)|}{\|f\|_{[-1,1]}}}
\leq c_2 \min\{B_{n,k,x}^*, M_{n,k}^*\}$$
for all $x \in (-1,1)$, where $c_1 > 0$ and $c_2 > 0$ are absolute constants.
This result should be compared with the inequalities
$$c_1 \min\{B_{n,k,x}, M_{n,k}\}
\leq \sup_{f \in {\Cal P}_{n,k}}{\frac {|f^\prime(x)|}{\|f\|_{[-1,1]}}}
\leq c_2 \min\{B_{n,k,x}, M_{n,k}\}$$
for all $x \in (-1,1)$, where $c_1 > 0$ and $c_2 > 0$ are absolute constants.
See [4] and [11].
It may be surprising that there is a significant difference between the real and
complex cases as far as Markov- and Bernstein-type inequalities are concerned.
In [3] essentially sharp Markov- and Bernstein-type inequalities for the classes 
${\Cal P}_{n,k}$ are proved even in $L_p$ norms on $[-1,1]$ for all $p > 0$.   

In this paper we revisit Erd\H os's paper [16] and make several remarks to 
his Markovi-type inequality in it. Erd\H os claimed in [16] that 
his method gave a Markov factor slightly better than $en/2$, namely, 
$$\|f^{\prime}\|_{[-1,1]} \leq \frac n2 \left( \frac{n}{n-1} \right)^{n-1} \|f\|_{[-1,1]}$$
for all $f \in {\Cal P}_n$ having all their zeros in ${\Bbb R} \setminus (-1,1)$.
Indeed, at some points of his arguments, by replacing applications of the inequality 
$1+x \leq e^x$ with an application of the inequality between the geometric and 
arithmetic means of nonnegative numbers, we can easily see this slight improvement.    

In 1963 Lorentz [20] proved that there is an absolute constant $c > 0$ such that 
$$|f^\prime(x)| \leq 
c \, \min \left\{ \sqrt{\frac{n}{1-x^2}}\,, n \right\} |f|_{[-1,1]}\,, \qquad x \in (-1,1)\,,$$
for all $f \in {\Cal B}_n(-1,1)$, where
$${\Cal B}_d(a,b) := \left\{ f: f(x) = \sum_{j=0}^d{a_j(b-x)^j(x-a)^{d-j}}\,, \quad a_j \geq 0, \quad j=0,1,\ldots, d \right\}\,.$$
for real numbers $a \leq b$ and nonnegative integers $d$. 
He also made the observation that if $f \in{\Cal P}_{n,0}$ then either $f \in {\Cal B}_n(-1,1)$ 
or $-f \in {\Cal B}_n(-1,1)$, 
where ${\Cal P}_{n,0}$ denotes the collection of all $f \in {\Cal P}_n$ having all 
their zeros outside the open unit disk. Scheick [23] has found the best possible constant $c$ in Lorentz's Markov-type inequality for $f \in {\Cal B}_n(-1,1)$. He showed that 
$$\|f^\prime\|_{[-1,1]} \leq \frac{en}{2} \|f\|_{[-1,1]}$$
for all $f \in {\Cal B}_n(-1,1)$.

An elementary, but very useful tool for proving inequalities for polynomials with 
restricted zeros is the {\it Bernstein} or {\it Lorentz representation} of polynomials.
Namely, as Lorentz observed it, if $p \in \Cal{P}_{n,0}$ is positive on $(-1,1)$ then 
it is of the form 
$$f(x) = \sum_{j=0}^d{a_j(1-x)^j(x+1)^{d-j}}\,, \qquad a_j \geq 0\,, \quad j=0,1,\, \ldots , d \,, \tag 1.1$$ 
with $d=n$. This is formulated  as Lemma 3.1 in this paper and its simple proof is reproduced. 
Moreover, if a polynomial $p \in \Cal{P}_n$ is positive on $(-1,1)$ and has no zeros in the 
ellipse  $L_\varepsilon$ with large axis $[-1,1]$ and small axis 
$[-\varepsilon i, \varepsilon i]$ ($\varepsilon \in [-1,1]$) then it has a Lorentz 
representation (1.1) with $d \leq 3n\varepsilon^{-2}$.
See [14]. Combining this with Lorentz's Markov- and Bernstein-type inequality gives that  
there is an absolute constant $c > 0$ such that
$$|p^{\prime}(x)| \leq c \,\min \left\{\frac{\sqrt{n}}{\varepsilon \sqrt{1-x^2}}\,,
\enskip \frac{n}{\varepsilon^2} \right\} \,\|p\|_{[-1,1]}\,, \qquad x \in (-1,1)\,,$$
for all $p \in {\Cal{P}_n}$ having no zeros in $L_\varepsilon$.

The minimal value of $d \in {\Bbb N}$ for which a polynomial $f$ has a representation (1.1) is called
the {\it Lorentz degree} of the polynomial and it is denoted by $d(f)$.
It follows from the already mentioned result in [14] that $d(p) < \infty$ if and only if
$p$ has no zeros in $(-1,1)$. This is a theorem ascribed to Hausdorff.
In addition, it has been proved in [8] that if
$$p(x) = ((x-a)^2 + \varepsilon^2(1-a^2))^n, \qquad 0 < \varepsilon \leq 1, \quad -1 < a < 1\,,$$
then
$$c_1n \varepsilon^{-2} \leq d(p) \leq c_2n \varepsilon^{-2}$$
with absolute constants $c_1 > 0$ and $c_2 > 0$.
Lorentz degree of trigonometric polynomials on an interval $(-\omega,\omega)$ shorter than 
the period is studied in [15].

\head 2. New Results \endhead

For $p > 0$ let
$$\|f\|_p := \bigg( \int_{-1}^1{|f(x| \, dx} \bigg)^{1/p}\,, \qquad \|f\|_{\infty} := \max_{x \in [-1,1]}{|f(x)|} \,.$$
As in Section 1 we will use the following notation. Let ${\Cal P}_n$ denote the 
collection of all polynomials of degree at most $n$ with {\it real coefficients}.
For real numbers $a \leq b$ and $d \in {\Bbb N}$ let
$${\Cal B}_d(a,b) := \left\{ f: f(x) = \sum_{j=0}^d{a_j(b-x)^j(x-a)^{d-j}}\,, \quad a_j \geq 0, \quad j=0,1,\ldots, d \right\}\,.$$
Let ${\Cal P}_{n,0}$ denote the collection of all $f \in {\Cal P}_n$ 
having all their zeros outside the open unit disk. Our first two results are the 
right Nikolskii-type inequalities for the classes ${\Cal B}_d(-1,1)$ and ${\Cal P}_{n,0}$.

\proclaim{Theorem 2.1} We have
$$\|f\|_p \leq \left( \frac{qd+1}{2} \right)^{1/q-1/p} \|f\|_q$$
for all $f \in {\Cal B}_d(-1,1)$ and for all $0 < q < p \leq \infty$.
Equality holds only for $f(x) := c(x \pm 1)^n$ with a constant $c \geq 0$.
\endproclaim

Combining Theorem 2.1 with Lemma 3.2 gives the following.

\proclaim{Theorem 2.2} We have
$$\|f\|_p \leq \left( \frac{qn+1}{2} \right)^{1/q-1/p} \|f\|_q$$
for all $f \in {\Cal P}_{n,0}$ and for all $0 < q < p \leq \infty$.
Equality holds only for $f(x) := c(x \pm 1)^n$ with a constant $0 \neq c \in {\Bbb R}$.
\endproclaim

An application of Theorem 2.1 with $q=1$ and $p=\infty$ allows us to prove the
following a sharp Markov-type inequality for all $f \in {\Cal P}_n$ such that
$f^\prime \in {\Cal B}_{d-1}(-1,1)$.

\proclaim{Theorem 2.3} We have
$$\|f^\prime\|_{\infty} \leq d \, \|f\|_{\infty}$$
for all $f \in {\Cal P}_d$ for which $f^\prime \in {\Cal B}_{d-1}(-1,1)$.
Equality holds only for $f(x) := c(x \pm 1)^n$ with a constant $c \geq 0$.
\endproclaim

Combining Theorem 2.3 with Lemma 3.2 gives the following.

\proclaim{Theorem 2.4} We have
$$\|f^\prime\|_{\infty} \leq n \, \|f\|_{\infty}$$
for all $f \in {\Cal P}_n$ for which $f^\prime$ has all its zeros outside
the open unit disk.
Equality holds only for $f(x) := c((x \pm 1)^n - 2^{n-1})$ with a constant $0 \neq c \in {\Bbb R}$.
\endproclaim

Our final result is a sharp Markov-type inequality for all $f \in {\Cal P}_n$ 
which are monotone on $[-1,1]$ and have all their zeros in ${\Bbb R} \setminus (-1,1)$. 
Erd\H os claimed this in [16] but he did not give a hint how to prove this. 
Experts seem to be  puzzled by this observation of Erd\H os even today.

\proclaim{Theorem 2.5} We have
$$\|f^\prime\|_{\infty} \leq \frac n2 \, \|f\|_{\infty}$$
for all $f \in {\Cal P}_n$ which is monotone on $[-1,1]$ and has all its zeros
in ${\Bbb R} \setminus (-1,1)$.
Equality holds only for $f(x) := c(x \pm 1)^n$ with a constant $0 \neq c \in {\Bbb R}$.
\endproclaim

We note that there is a incorrect hint to Part c] if Exercise 10
on page 482 of the book [2] suggesting that Theorem 2.5 holds.
However, it was discovered by M. Boedihardjo that the hint to part c] of E.10 
on page 482 of the book [2] does not work out. Here we claim a proof of 
Theorem 2.5 as a consequence of Theorem 2.4. A direct elementary proof of Theorem 2.5 
by using undergraduate calculus would be desirable.

\head 3. Lemmas \endhead

\proclaim{Lemma 3.1}
Let $a \leq c \leq d \leq b$ be real numbers, and let $d$ be a nonnegative integer. 
Then ${\Cal B}_d(a,b) \subset {\Cal B}_d(c,d)$.
\endproclaim  

\demo{Proof of Lemma 3.1}
This follows from the identities
$$x-a = \frac{c-a}{d-c} \, (x-c) + \frac{d-a}{d-c} \, (d-x)$$
and 
$$b-x = \frac{b-c}{d-c} \, (x-c) + \frac{b-d}{d-c} \, (d-x)$$
valid for all $x \in {\Bbb C}$.
\qed \enddemo

\proclaim{Lemma 3.2}
Suppose $f \in {\Cal P}_n$ has all its zeros outside the open unit disk. Then either 
$f \in {\Cal B}_n(-1,1)$ or $-f \in {\Cal B}_n(-1,1)$. 
\endproclaim

\demo{Proof of Lemma 3.2}
This follows from the identities
$$x-\alpha = \frac{1-\alpha}{2} \, (x+1) - \frac{\alpha+1}{2} \, (1-x)$$
and 
$$(x-\alpha)(x-\overline{\alpha}) = \frac 12 |1+\alpha|^2 (1-x)^2 + \frac 12 (|\alpha|^2-1)(1-x^2) + \frac 12 |1-\alpha|^2 (x+1)^2$$
valid for all $x \in {\Bbb C}$ and $\alpha \in {\Bbb C}$.
\qed \enddemo

\proclaim{Lemma 3.3} We have 
$$(\max\{f(a),f(b)\})^q \leq \frac{qd+1}{b-a} \int_a^b{f(x)^q \, dx}$$
for all $f \in {\Cal B}_d(a,b)$.
\endproclaim

\demo{Proof of Lemma 3.3} 
Let $f \in {\Cal B}_d(a,b)$ be of the form 
$$f(x) = \sum_{j=0}^d{a_j(b-x)^j(x-a)^{d-j}}\,, \qquad a_j \geq 0, \quad j=0,1, \ldots ,d\,.$$ 
Then
$$\split f(b)^q = & \, (a_0(b-a))^{dq} = \frac{qd+1}{b-a} \, \int_a^b{(a_0(b-x))^{dq} \, dx} \cr
\leq & \, \frac{qd+1}{b-a} \, \int_a^b{\Bigg( \sum_{j=0}^d{a_j(b-x)^j(x-a)^{d-j}} \Bigg)^q\, dx} \cr 
\leq & \, \frac{qd+1}{b-a} \int_a^b{f(x)^q \, dx}\,. \cr \endsplit $$
Similarly, 
$$\split f(a)^q = & \, (a_d(b-a))^{dq} = \frac{qd+1}{b-a} \, \int_a^b{(a_d(x-a))^{dq} \, dx} \cr
\leq & \, \frac{qd+1}{b-a} \, \int_a^b{\Bigg( \sum_{j=0}^d{a_j(b-x)^j(x-a)^{d-j}} \Bigg)^q\, dx} \cr 
\leq & \, \frac{qd+1}{b-a} \int_a^b{f(x)^q \, dx}\,. \cr \endsplit$$
\qed \enddemo

\proclaim{Lemma 3.4} We have
$$\|f\|_{\infty}^q \leq \frac{qd+1}{2} \, \|f\|_q^q\,.$$
for all $f \in {\Cal B}_d(-1,1)$.
\endproclaim

\demo{Proof of Lemma 3.4}
Let $y \in [-1,1]$ be such that $f(y) = \|f\|_{\infty}$. 
By Lemma 3.1 we have 
$${\Cal B}_d(-1,1) \subset {\Cal B}_d(-1,y) \cap {\Cal B}_d(y,1)\,.$$ 
Hence Lemma 3.2 yields
$$(y+1)f(y)^q \leq (qd + 1) \int_{-1}^y{f(x)^q \, dx}$$
and
$$(1-y)f(y)^q \leq (qd + 1) \int_{y}^1{f(x)^q \, dx}$$
Adding the above two inequalities, we conclude
$$\|f\|_{\infty}^q = f(y)^q \leq \frac{qd + 1}{2} \, \int_{-1}^1{f(x)^q \, dx} = \frac{qd + 1}{2} \, \|f\|_q^q\,.$$
\qed \enddemo

\head Proof of the Theorems \endhead

\demo{Proof of Theorem 2.1}
When $p = \infty$ the Theorem follows from Lemma 3.3. Now let 
Let $f \in {\Cal B}_d(-1,1)$ and $0 < q \leq p \leq \infty$. Using Lemma 3.3 we obtain
$$\split \|f\|_p^p = & \, \int_{-1}^1{f(x)^p \, dx} \leq \Bigg( \int_{-1}^1{f(x)^q \, dx} \Bigg) \|f\|_{\infty}^{p-q}
\leq \, \|f\|_q^q \left( \frac{qd+1}{2} \right)^{(p-q)/q} \|f\|_q^{p-q} \cr
= & \, \left( \frac{qd+1}{2} \right)^{(p-q)/q} \|f\|_q^p \,, \cr \endsplit $$
hence
$$\|f\|_p \leq \left( \frac{qd+1}{2} \right)^{1/q - 1/p} \|f\|_q \,.$$
\qed \enddemo

\demo{Proof of Theorem 2.2}.
Combining Theorem 2.1 and Lemma 3.2 gives the result.
\qed \enddemo

\demo{Proof of Theorem 2.3}
Applying Theorem 2.1 with $f$ replaced by $f^\prime \in {\Cal B}_{d-1}$, $p := \infty$ and $q := 1$, we obtain 
$$\|f^\prime\|_{\infty} \leq \frac d2 \, \int_{-1}^1{f^\prime(x) \, dx} = 
\frac d2 \, (f(1) - f(-1)) \leq d \, \|f\|_{\infty} \,.$$ 
\qed \enddemo

\demo{Proof of Theorem 2.4}
Assume that $f^\prime \in {\Cal P}_{n-1}$ has no zeros in the open unit disk. Then, 
by Lemma 3.2 either $f^\prime \in {\Cal B}_{n-1}(-1,1)$ or $-f^\prime \in {\Cal B}_{n-1}(-1,1)$. Without loss of generality we may assume that $f^\prime \in {\Cal B}_n(-1,1)$, 
and Theorem 2.3 gives the result. 
\qed \enddemo

\demo{Proof of Theorem 2.5}
Assume that $f \in {\Cal P}_n$ is monotone on $[-1,1]$ and has all its zeros in 
${\Bbb R} \setminus (-1,1)$. Then, by Rolle's Theorem $f^\prime$ has all its zeros in 
${\Bbb R} \setminus (-1,1)$, and hence Lemma 3.2 implies that  either 
$f^\prime \in {\Cal B}_{n-1}(-1,1)$ or $-f^\prime \in {\Cal B}_{n-1}(-1,1)$. 
Without loss of generality we may assume that $f^\prime \in {\Cal B}_n(-1,1)$.
Applying Theorem 2.3 we conclude that 
$$\|f^\prime\|_{\infty} \leq \frac n2 \, \int_{-1}^1{f^\prime(x) \, dx} = 
\frac n2 \, (f(1) - f(-1)) \leq \frac n2 \, \|f\|_{[-1,1]} \,,$$     
where in the last step we used that $f(1)f(-1) \geq 0$ since $f \in {\Cal P}_n$ has all its zeros in ${\Bbb R} \setminus (-1,1)$.
\qed \enddemo

\enddocument